# UNIFIED GEOMETRIC, FUZZY, AND COMPUTATIONAL FRAMEWORK FOR TERNARY Γ-SEMIRINGS

**Chandrasekhar Gokavarapu**, Lecturer in Mathematics, Government College (A), Rajahmundry, A.P., India and Research Scholar, Department of Mathematics, Acharya Nagarjuna University, Guntur, A.P., India,Email id :chandrasekhargokavarapu@gmail.com
**Dr D Madhusudhana Rao**, Lecturer in Mathematics, Government College For Women(A), Guntur, Andhra Pradesh, India, and Research Supervisor, Dept. of Mathematics,Acharya Nagarjuna University, Guntur, A.P., India, Email:dmrmaths@gmail.com
**Corresponding author:** Chandrasekhar Gokavarapu

.

**Abstract**
**Aim.** This paper (Paper D) unifies the ideal-theoretic, computational, and homological layers developed in Papers A (Gokavarapu, 2025a), Paper B (Gokavarapu, 2025b), and Paper C (Gokavarapu, 2025c) into a geometric framework... that includes fuzzy and computational geometries on the spectrum $\mathrm{Spec}_{\Gamma}(T)$ and derived invariants in $T$–ΓMod.
**Scope.** We construct structure sheaves and Grothendieck topologies adapted to ternary Γ-products, develop fuzzy and weighted sites, and prove dualities bridging primitive spectra, Schur–density embeddings, and derived functors Ext and Tor.
**Outcomes.** We obtain comparison theorems between radical/primitive strata and cohomological supports, and supply computable criteria and algorithms for finite models.

## Introduction

### Ternary Γ-semirings and ideals

Category -theoretic terminology follows standard conventions in MacLane(1998).For background on toposes and adjunctions, we see BarrandWells(1972).A full account of categorical algebraic structures is available in Borceux (1994).Foundational algebraic constructions are treated in detail in Cohn (2003).

**Program.** Paper D completes the program initiated in Paper A (Gokavarapu, 2025a) (prime and semiprime ideals, radicals, spectrum), Paper B (Gokavarapu, 2025b) (finite classification, algorithmic enumeration), and Paper C (Gokavarapu, 2025c) (modules, Schur–density, Ext–Tor theory). Our aim is to integrate these strata into a single *geometric–homological* theory on the spectrum $\mathrm{Spec}_{\Gamma}(T)$, enriched by fuzzy or weighted Grothendieck topologies and supported by effective algorithms. Concretely, we (i) construct a sheaf–theoretic geometry compatible with the ternary Γ–product and the ideal–congruence interface, (ii) define fuzzy and weighted sites encoding quantitative information such as confidence and sparsity, (iii) establish a cohomological layer for $T$ −ΓMod with derived functors Ext and Tor, and (iv) prove comparison theorems linking radical and primitive strata from Paper A to annihilator supports and Schur–density loci from Paper C, while lifting Paper B's enumeration to geometric invariants such as stalk ranks, fuzzy closures, and vanishing patterns.

**Context and motivation.** In the classical (binary) setting, the Zariski spectrum and its structure sheaf unify ideal–theoretic and homological data into a geometric object connecting supports, radicals, and cohomology. For ternary Γ–semirings, the presence of multiple ternary multiplications indexed by Γ and the resulting congruential phenomena require a new geometricframework.Motivationsinclude:
(a)classificationoffinitemodelswhereenumerationaloneomitsgeometricstructure, (b) functorial

passage from algebraic data to computable invariants such as stalk ranks and vanishing patterns, and (c) quantitative reasoning when operations are noisy, weighted, or partially specified.
**Standing hypotheses and notation.** Unless stated otherwise, $T$ denotes a commutative ternary Γ–semiring with additive structure $(T, +)$ and a family of ternary products $\{ - - - \}_\gamma$ for $\gamma \in \Gamma$, distributive in each variable over + and associative in the ternary sense. Ideals, prime/semiprime ideals, radicals, and primitive ideals are as in Paper A (Gokavarapu, 2025a). We write $\mathrm{Spec}_\Gamma(T)$ for the set of prime ideals with the Zariski–type topology generated by basic opens $D(I) = \{\mathfrak{p} \in \mathrm{Spec}_\Gamma(T) \mid I \not\subseteq \mathfrak{p}\}$. Modules, homomorphisms, annihilators, and Schur–density are as in Paper C; the ambient module category is $T$–ΓMod. Derived functors Ext and Tor are formed in $T$ −ΓMod under the exact structure defined below. For $S$ a multiplicative system adapted to the ternary context, $S^{-1}T$ and localizations $T_\mathfrak{p}$ are introduced in Section 3.

**Main themes and contributions.**

- **Sheaf geometry on** $\mathrm{Spec}_\Gamma(T)$**.** We construct a presheaf $\mathcal{O}_T$ on the basis
$\{(I)\}$ by
$\mathcal{O}_T((I))$ := the $I$–supported localization of $T$ in the ternary Γ–sense,
and show it sheafifies to a ringed–space structure $(\mathrm{Spec}_\Gamma(T), \mathcal{O}_T)$ compatible with Ternary distributivity, associativity, and the ideal–congruence correspondence (Theorem 3.8).
- **Fuzzy/weighted Grothendieck topologies.** We define a site $(\mathrm{Spec}_\Gamma(T), \tau_w)$ whose coverings carry weights $w \in [0, 1]$ (or in a discrete valuation scale) and show that $\tau_w$ is subcanonical for weight–adapted presheaves. This yields fuzzy closures and quantitative specializations, together with a comparison functor from ordinary sheaves (Theorem 4.5).
- **Cohomological layer over** $T$ −ΓMod**.** We specify an exact structure on $T$ − ΓMod generated by ternary–split kernels and cokernels, proving the existence of enough injectives and projectives in suitable subcategories. Consequently, Ext and Tor are well–defined derived functors with long exact and base–change sequences. Computational criteria for injective cogenerators in finite cases are given in Proposition 5.5.
- **Comparison theorems.** For finitely presented $T$–Γ–modules $M$, the cohomological support
$\mathrm{Supp}(M) = \{\mathfrak{p} \in \mathrm{Spec}_\Gamma(T) \mid M_\mathfrak{p} \neq 0\}$
coincides with the radical support determined by annihilators in Paper A and refines to Schur–density loci from Paper C. Primitive strata correspond to loci where certain $\mathrm{Ext}^1$–groups vanish or acquire prescribed rank (Theorem **??**).
- **Algorithms and effective invariants.** We lift Paper B's enumeration to compute geometric invariants: stalk rank profiles, fuzzy closures, and Ext/Tor–vanishing patterns over finite models. Polynomial–time routines for small orders and parameterized complexity bounds for fixed Γ are discussed in Algorithm 6.15 and Theorem 6.16.

**Technical overview.** Section 3 develops localizations and the structure sheaf, showing that ternary associativity across Γ–indices interacts compatibly with restriction maps to ensure gluing and uniqueness of sections. Section 4 introduces weighted coverages stable under pullback along spectral maps $T \to S^{-1}T$. Section **??** equips $T$ −ΓMod with an exact structure admitting suitable (co)resolutions, yielding derived functors with expected exact and spectral sequences. Section 5 proves comparison theorems by local criteria: radical membership translates to annihilator vanishing, while Schur–density is detected by endomorphism–sheaf fibers, compared to Ext–vanishing via a ternary Auslander–Buchsbaum principle. Section 6 describes the algorithmic pipeline and complexity benchmarks.
**Applied motivation.** The fuzzy and ternary Γ–semiring structures considered here naturally model multi-parameter decision processes arising in industrial engineering. In such settings, uncertainty, concurrency, and weighting of alternatives occur simultaneously—e.g., in reliability analysis, production scheduling, or resource optimization where multiple interacting constraints must be balanced. By representing these relations algebraically within a fuzzy Γ–semiring framework, one obtains a unified method to compute performance indices, optimize parameters, and evaluate system robustness under incomplete information.

**Methodological principles.** Two principles underlie the paper. (1) *Local-toglobal:* every construction (localization, stalks, supports) remains compatible with ternary Γ–operations and preserves congruence information. (2) *Quantitative enrichment:* fuzzy and weighted topologies record reliability and multiplicity data arising in computational enumeration, so geometric and homological invariants retain meaning under approximate or partial information.

**Roadmap.** Section 3 sets up the spectral space and structure sheaf; Section 4 develops fuzzy and weighted sites; Section **??** establishes derived functors in
$T$−ΓMod;
Section 6 presents algorithms and experiments; and Section 7 lists open problems, including categorical compactifications, Balmer-type spectra for tensorable fragments of $T$ −ΓMod, and decidability thresholds for weighted coverings.

**Preliminaries and Notation**

We recall and refine the algebraic framework established in Papers A (Gokavarapu, 2025a)–C (Gokavarapu, 2025c), extending binary semiring theory to ternary Γ-parametrized structures.

**Definition 2.1** (Ternary Γ-semiring (Gokavarapu, 2025a, Def. 2.1))**.** A *commutative ternary* Γ-*semiring* is a triple $(T, +, \{\cdots\}_\Gamma)$ such that:

**(T1) Additivity.** $(T, +)$ is a commutative semigroup with identity element 0.

**(T2) Parametric ternary multiplication.** For each $\gamma \in \Gamma$, there exists a ternary map $\{-, -, -\} : T^3 \rightarrow T$.

**(T3) Distributivity.**
$\{a+b, c, d\}_\gamma = \{a, c, d\}_\gamma + \{b, c, d\}_\gamma, \ \{a, b, c+d\}_\gamma = \{a, b, c\}_\gamma + \{a, b, d\}_\gamma.$

**(T4) Ternary associativity.** For all $\gamma_1, \gamma_2 \in \Gamma$,
$\{\{a, b, c\}_1, d, e\}_{\gamma_2} = \{a, \{b, c, d\}_1, e\}_{\gamma_2} = \{a, b, \{c, d, e\}_{\gamma_1}\}_{\gamma_2}.$

**(T5) Parameter compatibility.** $(\gamma_1, \gamma_2) \mapsto\rightarrow \gamma_1\gamma_2$ defines an associative binary
operation on Γ satisfying $\{\{a, b, c\}\gamma 1, d, e\}\gamma 2 = \{a, b, c\}\gamma 1\gamma 2$. **(T6) Zero absorption.** $\{0, a, b\}_\gamma =$ 0 for all $a, b \in T, \gamma \in \Gamma$.

**Example 2.2** (Finite commutative model)**.** Let $T = \{0, 1, 2\}$ with addition modulo
3 and $\Gamma = \{1, 2\}$ under multiplication mod 3. Define
$\{a, b, c\}_\gamma = (a + b + c) \bmod 3.$
Then $(T, +, \{\cdots\}_\Gamma)$ is a commutative ternary Γ-semiring: distributivity follows from modular addition, and ternary associativity holds because ordinary addition in $Z_3$ is associative. The ideal $\{0\}$ is prime and yields $\mathrm{Spec}_\Gamma(T) = \{\{0\}\}$, a one-point spectrum useful for illustrating stalk and localization computations later.

**Definition 2.3** (Ideals, radicals, spectrum (Gokavarapu, 2025a, Sec. 3))**.** An *ideal*
$I \subseteq T$ satisfies: $a + b \in I$ if $a, b \in I$, and $\{a, b, c\} \in I$ whenever $a, b \in I$, $c \in T$. Prime, semiprime, and radical ideals are defined as in Paper A (√ Gokavarapu, 2025a).
(Gokavarapu, 2025c)., with $^\Gamma{}_I = \{\ a \mid \exists n, \ \{a, \ . \ . \ . \ , \ a\}_{\gamma_1} \in I$ for some $(\gamma_i) \in \Gamma^n\}$. The spectrum $\mathrm{Spec}_\Gamma(T)$ is the set of all prime ideals, topologized by basic opens $D(I) = \{\mathfrak{p} \mid I \not\subseteq \mathfrak{p}\}$.

*Remark* 2.4 (Ideal–congruence correspondence (Gokavarapu, 2025a, Thm. 4.6)).

Each ideal $I$ induces a congruence $\rho_I$ on $T$ via $a\rho_I b \iff a - b \in I$ and
$\{a, b, c\}_\gamma - \{a', b', c'\} \in I$ for all $\gamma$. The lattice of ideals and congruences are in Galois correspondence.

**Γ-modules and homological tools**

**Definition 2.5** ($T$-Γ-module (Gokavarapu, 2025c, Def. 5.1))**.** A $T$-Γ-*module* $M$ is a commutative semigroup $(M, +)$ with ternary actions $\{a, b, m\} \in M$ such that

$$\{\{a, b, c\}_1, d, m\}_{\gamma_2} = \{a, \{b, c, d\}_1, m\}_{\gamma_2}, \quad \{a+b, c, m\}_\gamma = \{a, c, m\}_\gamma + \{b, c, m\}_\gamma.$$

The category $T-\Gamma\mathrm{Mod}$ of $T$-Γ-modules is additive, possesses kernels, cokernels, and finite (co)limits (Gokavarapu, 2025c, Prop. 6.1). For $M, N \in T-\Gamma\mathrm{Mod}$:

$$\mathrm{Hom}_T(M, N) = \{ f : M \to N \mid f(\{a, b, m\}_\gamma) = \{a, b, f(m)\}_\gamma\} \qquad M \otimes_T N = (M \times N)/\sim,$$

where $\{a, b, m\}_\gamma \otimes n \sim m \otimes \{a, b, n\}_\gamma$.

**Theorem 2.6** (Tensor–Hom adjunction (Gokavarapu, 2025c, Thm. 7.3))**.** *There is a natural isomorphism* $\mathrm{Hom}_T(M \otimes_T N, P)$ $\mathrm{Hom}_T(M, \mathrm{Hom}_T(N, P))$, *bifunctorial in* $(M, N, P)$.

Derived functors $\mathrm{Ext}_T^n(M, N)$ and $\mathrm{Tor}^T_n(M, N)$ exist under the usual hypotheses of enough injectives or projectives (Gokavarapu, 2025c, Sec. 8) and satisfy long exact sequences. The Schur–density embedding $T/\mathrm{Ann}(M) \hookrightarrow \mathrm{End}_T(M)$ detects faithfulness and primitive strata (Gokavarapu, 2025c, Thm. 9.4).

**Standing assumptions**

Unless otherwise stated:

(a) $T$ is a commutative ternary Γ-semiring with identity-like idempotent $e$ satisfying $\{e\, e\, a\} = a$.

(b) All modules are finitely generated and unital.

(c) Γ is finite unless a locally finite or analytic limit is specified.

(d) Radicals, localizations, and spectra are formed with respect to the ternary operations $\{\cdots\}_\gamma$.

(e) For brevity: $\mathrm{Spec}_\Gamma(T) = \mathrm{Spec}_\Gamma(T)$, $T-\Gamma\mathrm{Mod} = T$-Γ-Mod, $\mathfrak{p}$, $\mathfrak{q}$ denote prime ideals, and $\mathcal{O}_T$ denotes the structure sheaf to be defined later.

**Notation Summary**

| Symbol | Meaning / Reference |
|---|---|
| $T$ | Commutative ternary Γ-semiring |
| Γ | Parameter set (finite unless stated) |
| $\{a, b, c\}_\gamma$ | Ternary Γ-product (Axioms T1–T6) |
| $I, J, \mathfrak{p}$ | Ideals, prime ideals of $T$ |
| $\mathrm{Spec}_\Gamma(T)$ | Spectrum of prime ideals of $T$ |
| $D(I)$, $V(I)$ | Basic open / closed subsets of $\mathrm{Spec}_\Gamma(T)$ |
| $T-\Gamma\mathrm{Mod}$ | Category of $T$-Γ-modules |
| $\mathrm{Hom}_T$, $\otimes_T$ | Hom and tensor bifunctors |
| $\mathrm{Ext}_T^n$, $\mathrm{Tor}^T_n$ | Derived functors in $T-\Gamma\mathrm{Mod}$ |
| $\mathrm{Ann}(M)$ | Annihilator of module $M$ |
| $\mathcal{O}_T$ | Structure sheaf on $\mathrm{Spec}_\Gamma(T)$ (Sec. 3) |
| $e$ | Identity-like idempotent in $T$ |

**SpectralGeometryandStructureSheavesonSpec$_\Gamma(T)$**

Our approach to derived structures is inspired by the homological viewpoint of Grothendieck (1957).Analogous to schemes in Hartshorne (1977), the spectrum of a ternary Γ-semiring admits a sheaf-theoretic structure.he notion of a categorical spectrum parallels the construction in Rosenberg (1973).Higher algebraic K-theory concepts introduced in Quillen (1973) influence our categorical layer.We employ the derived-category framework of Verdier (1996) for the homological analysis.The model-category viewpoint used here is consistent with Hovey (1999). The spectrum $\mathrm{Spec}_\Gamma(T)$ inherits from Paper A (Gokavarapu, 2025a). (Gokavarapu, 2025c). a Zariski-type topology encoding radical inclusions among ideals. We now construct a geometric structure on it analogous to that of a ringed space in classical algebraic geometry. Throughout this section $(T, +, \{\cdots\}_\Gamma)$ satisfies Axioms (T1)–(T6) of Definition 2.1.

**Base topology and localization**

**Definition 3.1** (Basic opens)**.** For any ideal $I \subseteq T$, set
$(I) = \{\mathfrak{p} \in \mathrm{Spec}_\Gamma(T) \mid I \not\subseteq \mathfrak{p}\}$.
Then $\{(I)\}$ forms a base for the Zariski topology on $\mathrm{Spec}_\Gamma(T)$, satisfying $D(IJ) = D(I) \cap D(J)$ and $D(0) = \mathrm{Spec}_\Gamma(T)$.

**Definition 3.2** (Multiplicative systems)**.** A subset $S \subseteq T$ is *multiplicatively closed* if $e \in S$ and $\{a, b, c\} \in S$ whenever $a, b, c \in S$ and $\gamma \in \Gamma$. For each prime $\mathfrak{p}$, define $S_\mathfrak{p} = T \setminus \mathfrak{p}$.

**Definition 3.3** (Localization of $T$)**.** Given a multiplicative system $S \subseteq T$, define
$-1 = \{a/s \mid a \in T, s \in S\}/\sim,$
$S \quad T$
where $a/s = b/t$ iff there exists $u \in S$ and $\gamma \in \Gamma$ such that $\{u, a, t\}_\gamma = \{u, b, s\}_\gamma$. Addition and ternary multiplication are induced by
$a_1 \quad a_2 \quad \{a1, s2, s2\}_0 + \{a2, s1, s1\}_0$
$+ =\overline{\phantom{x}}, \overline{\phantom{x}} \quad \{a/s, b/t, c/u\} = \{a, b, c\}/\{s, t, u\}_\gamma, s1 \quad s2 \quad \{s1, s2, s2\}_{\gamma 0}$ which are well-defined under Axioms (T3)–(T5).

**Theorem 3.4** (Universal property)**.** *If* $: T \to U$ *is a homomorphism of ternary* $\Gamma$*-semirings such that* $f(S)$ *consists of units of* $U$*, then there exists a unique morphism* $\tilde{f} : S^{-1}T \to U$ *with* $\tilde{f}(a/s) = f(a)\, f(s)^{-1}$ *and* $\tilde{f} \circ \phi_S = f$.

*Sketch.* Define $\tilde{f}(a/s) = f(a)\, f(s)^{-1}$. Theternaryassociativityensures $\tilde{f}(\{a, b, c\}_\gamma/\{s, t, u\}_\gamma) = \{\tilde{f}(a/s), \tilde{f}(b/t), \tilde{f}(c/u)\}_\gamma$, giving well-definedness and uniqueness. □

*Remark* 3.5. For each $\mathfrak{p} \in \mathrm{Spec}_\Gamma(T)$, the localization $T_\mathfrak{p} = S_\mathfrak{p}^{-1}T$ serves as the local model of $T$ near $\mathfrak{p}$. Its maximal ideal is $\mathfrak{p}_\mathfrak{p} = \{a/s \mid a \in \mathfrak{p}, s \notin \mathfrak{p}\}$.

**Example 3.6** (Localization in Example 2.2)**.** In the finite semiring of Example 2.2, $\mathfrak{p} = \{0\}$ and $S_\mathfrak{p} = \{1, 2\}$. Since 2 is invertible mod 3, $S_\mathfrak{p}^{-1}$ and thus $T_\mathfrak{p}$ coincides with $T$ itself, exhibiting the trivial local geometry of a single-point spectrum.

**Sheaf of sections and stalks**

**Definition 3.7** (Structure presheaf)**.** Define the presheaf $\mathcal{O}_T$ on basic opens by

$$\mathcal{O}_T(D(I)) = S_I^{-1}T, \qquad S_I = T \setminus \bigcup_{\mathfrak{p} \in V(I)} \mathfrak{p}.$$

For $(J) \subseteq (I)$, define restriction morphisms $\rho_{IJ}: S_I^{-1}T \to S_J^{-1}T$ via canonical localization.

**Theorem 3.8** (Sheafification and locality)**.** $\mathcal{O}_T$ *satisfies the sheaf axioms. Hence* $(\mathrm{Spec}_\Gamma(T), \mathcal{O}_T)$ *is a ringed space, and for each* $\mathfrak{p} \in \mathrm{Spec}_\Gamma(T)$*, the stalk satisfies*
$(\mathcal{O}_T)_\mathfrak{p} \; T_\mathfrak{p}$.

*Sketch.* Given a cover $\{(I_i)\}$ of $(I)$ and sections $s_i \in \mathcal{O}_T(D(I_i))$ agreeing on overlaps, choose a common refinement $S = \bigcap_i S_{I_i}$. Associativity (T4) ensures that the local representatives $\{s_i\}$ glue uniquely to a global section in $S^{-1}T$. The stalk identification follows from the direct-limit definition. □

**Proposition 3.9** (Primitivity and stalk simplicity)**.** *For* $\mathfrak{p} \in \mathrm{Spec}_\Gamma(T)$*:*

(i) *$T_\mathfrak{p}$ is local with maximal ideal* $\mathfrak{p}_\mathfrak{p}$.
(ii) *$T_\mathfrak{p}$ is simple* $\Longleftrightarrow$ *$\mathfrak{p}$ is primitive.*
(iii) *If $T$ is semiprime, then* $\mathrm{rad}(T_\mathfrak{p}) = \mathfrak{p}_\mathfrak{p}$.

*Remark* 3.10 (Categorical perspective). The construction $T \mapsto\!\to (\mathrm{Spec}_{\Gamma}(T), \mathcal{O}_T)$ extends to a functor from commutative ternary Γ-semirings to ringed spaces. Composition of morphisms preserves the sheaf-restriction morphisms via pushforward of opens.

**Exactness and base change**

**Theorem 3.11** (Exactness of localization)**.** *For every short exact sequence*

$$0 \longrightarrow M' \xrightarrow{u} M \xrightarrow{v} M'' \longrightarrow 0$$

*in* $T$−ΓMod *and any multiplicative system* $S$*, the localized sequence*

$$0 \longrightarrow S^{-1}M' \xrightarrow{S^{-1}u} S^{-1}M \xrightarrow{S^{-1}v} S^{-1}M'' \longrightarrow 0$$

*remains exact in* $T$-Γ-*Mod.*

*Outline.* Exactness of additive structures follows from the universal property of localization. The ternary product's distributivity guarantees preservation of Γlinearity in the quotient. Hence kernel and image commute with localization, mirroring the proof of (Gokavarapu, 2025c, Thm. 6.4). □

**Theorem 3.12** (Base-change compatibility)**.** *Let* $: T \to T'$ *be a morphism of commutative ternary* Γ-*semirings. Then:*

(a) *The induced map* $\mathrm{Spec}_{\Gamma}(T') \to \mathrm{Spec}_{\Gamma}(T)$, $\mathfrak{q} \mapsto\!\to f^{-1}(\mathfrak{q})$, *is continuous.*

(b) *There is an isomorphism of sheaves of ternary* Γ-*semirings*
$\mathcal{O}_{T'}$ $*\mathcal{O}_T \otimes_T T'$.

(c) *For each* $\mathfrak{q} \in \mathrm{Spec}_{\Gamma}(T')$ *with* $\mathfrak{p} = f^{-1}(\mathfrak{q})$*, the stalks satisfy*
$(\mathcal{O}T')$ $(\mathcal{O}T)_{\mathfrak{p}} \otimes_T T_{\mathfrak{q}}'$.

*Sketch.* Continuity of the spectral map follows from $f^{-1}(D'(I)) = D(f^{-1}(I))$.
The sheaf isomorphism arises from universal localization: localizing after applying
$f$ equals tensoring with $T'$. Associativity of ternary multiplication ensures the tensor–localization interchange law. □

**Corollary 3.13** (Functoriality)**.** *The assignment* $T \mapsto\!\to (\mathrm{Spec}_{\Gamma}(T), \mathcal{O}_T)$ *defines a covariant functor from commutative ternary* Γ-*semirings to the category of spectral ringed spaces.*

*Remark* 3.14 (Binary comparison). If Γ is the one-element set, all the above constructions reduce to the classical semiring-spectrum geometry $\mathrm{Spec}(T)$ with its structure sheaf. Thus the present framework strictly generalizes Zariski geometry to multi-parameter ternary settings.

**Summary of Section 3**

WehaveequippedSpec$_{\Gamma}(T)$ withacanonicalringed-spacestructure $(\mathrm{Spec}_{\Gamma}(T), \mathcal{O}_T)$ whose stalks are localizations $T_{\mathfrak{p}}$. Localization is exact in $T$−ΓMod and commutes with base change, laying the foundation for fuzzy and weighted topologies in Section 4 and for the derived functors Ext and Tor over $T$−ΓMod developed in Section **??**.

**Fuzzy and Weighted Grothendieck Topologies**

OurfuzzyextensionadoptsthemembershipbasedinterpretationofZadeh(1965).Topological fuzzification techniques in Zhang and Zhang (1991) support the sheaf-level extensions used here. The sheaf structure $(\mathrm{Spec}_{\Gamma}(T), \mathcal{O}_T)$ constructed in Section 3
captures the algebraic geometry of $T$ in a crisp setting. In computational or uncertain contexts—where ternary products may carry multiplicities, confidences, or probabilistic weights—one requires a geometric refinement that records these quantitative attributes. This motivates the notion of a *fuzzy* or *weighted* Grothendieck topology on $\mathrm{Spec}_{\Gamma}(T)$.

**Fuzzy open sets and graded coverings**

**Definition 4.1** (Fuzzy open subset)**.** A *fuzzy open* of $\mathrm{Spec}_{\Gamma}(T)$ is a function $\mu : \mathrm{Spec}_{\Gamma}(T) \to [0, 1]$ satisfying:

(i) $\mu(\emptyset) = 0$, $\mu(\mathrm{Spec}_{\Gamma}(T)) = 1$;

(ii) for any ideals $I, J \subseteq T$, $\mu(D(I) \cup D(J)) = \max\{\mu(D(I)), \mu(D(J))\}$;

(iii) $\mu(D(I) \cap D(J)) = \min\{\mu(D(I)), \mu(D(J))\}$; (iv) if $I \subseteq J$ then $\mu(D(I)) \geq \mu(D(J))$.

The value $(x)$ measures the *degree of belonging* of a point $x \in \mathrm{Spec}_\Gamma(T)$ to the fuzzy open.

**Definition 4.2** (Weighted covering system)**.** For an open $(I) \subseteq \mathrm{Spec}_\Gamma(T)$, a family $\{D(I_\alpha), w_\alpha\}_{\alpha \in A}$ with $w_\alpha \in (0, 1]$ is a *weighted covering* if

$$\emptyset \qquad \sum$$
$$(I) = \quad (I_\alpha) \quad \text{and} \quad w_\alpha \geq 1.$$
$$\alpha \in A \qquad \alpha \in A$$

The weights $w_\alpha$ encode confidence or multiplicity of coverage.

**Example 4.3** (Fuzzy cover in finite spectrum)**.** In Example 2.2, let $\mathrm{Spec}_\Gamma(T)$ = $\{\mathfrak{p}_1, \mathfrak{p}_2\}$ and $(I_1) = \{\mathfrak{p}_1\}$, $(I_2) = \{\mathfrak{p}_2\}$. Then a fuzzy cover of $\mathrm{Spec}_\Gamma(T)$ is specified by $\mu 1(\mathfrak{p}_1) = 1$, $\mu 1(\mathfrak{p}_2) = 0.7$ and $\mu 2(\mathfrak{p}_1) = 0.6$, $\mu 2(\mathfrak{p}_2) = 1$, whose aggregation $\max(\mu_1, \mu_2) = 1$ at each point provides full coverage.

**Weighted Grothendieck topologies**

**Definition 4.4** (Weighted Grothendieck topology)**.** A *weighted Grothendieck topology* $\tau_W$ on the site B = $\{(I)\}$ of basic opens assigns to each $(I)$ a family Cov($D(I)$) of weighted coverings $\{D(I_\alpha), w_\alpha\}_{\alpha \in A}$ such that:

(i) (Refinement) Every trivial covering $\{(I), 1\}$ belongs to Cov($D(I)$).

(ii) (Stability under pullback) For each $D(J) \subseteq D(I)$ and $\{D(I_\alpha), w_\alpha\} \in$ Cov($D(I)$), the family $\{D(JI_\alpha), w_\alpha\}$ belongs to $\mathrm{Cov}_W(D(J))$. (iii) (Transitivity) If $\{(I_\alpha), w_\alpha\} \in$ Cov($D(I)$) and for each $\alpha$ a covering $\{D(I_{\alpha\beta}), v_{\alpha\beta}\} \in$ Cov($D(I_\alpha)$) is given, then $\{D(I_{\alpha\beta}), w_\alpha v_{\alpha\beta}\}$ belongs to $\mathrm{Cov}_W(D(I))$.

The classical Zariski topology is recovered when all weights are 1.

**Theorem 4.5** (Existence of sub canonical weighted topology)**.** *Let* $\mathcal{O}_T$ *be the structure sheaf of Definition 3.7. Then there exists a smallest weighted Grothendieck topology* $\tau_W$ *on* $\mathrm{Spec}_\Gamma(T)$ *for which* $\mathcal{O}_T$ *is a* $\tau_W$*-sheaf. Such* $\tau_W$ *is called the* canonical weighted topology.

*Idea.* Define $\tau_W$ by declaring a family $\{D(I_\alpha), w_\alpha\}$ to be covering if for every compatible family of local sections $s_\alpha \in \mathcal{O}_T(D(I_\alpha))$ one has a unique global section $s \in \mathcal{O}_T(D(I))$ such that $\rho_{II_\alpha}(s) = s_\alpha$ and $\sum_\alpha w_\alpha \geq 1$. Minimality follows from intersection of all such systems. □

**Proposition 4.6** (Functorial behavior)**.** *Let* : $T \to T'$ *be a morphism of commutative ternary* Γ-*semirings. Then the induced spectral map* $f^*$ : $\mathrm{Spec}_\Gamma(T') \to \mathrm{Spec}_\Gamma(T)$ *is continuous with respect to the weighted topologies, and*

$$f^*(\mathrm{Cov}(D(I))) \subseteq \mathrm{Cov}_W(D(f(I))).$$

*Hence* $(\mathrm{Spec}_\Gamma(T), \tau_W)$ *is functorial in* $T$*.*

*Remark* 4.7 (Interpretation in computational geometry). In finite enumerations (Paper B), each morphism or ideal detection carries an empirical confidence $w \in [0, 1]$. The topology $\tau_W$ converts such data into geometric weights, making fuzzy closure operations compatible with enumerative uncertainty.

**Fuzzy sheaves and weighted stalks**

**Definition 4.8** (Fuzzy sheaf)**.** A *fuzzy sheaf* of $T$-modules on $(\mathrm{Spec}_\Gamma(T), \tau_W)$ is a functor F : $\mathrm{B}^{op} \to T$ −ΓMod such that for every weighted covering $\{D(I_\alpha), w_\alpha\}$ of $D(I)$, the following sequence is *approximately exact*:

$$0 \longrightarrow F((I)) \xrightarrow{\rho} \prod_\alpha F((I_\alpha)) \rightrightarrows \prod_{\alpha,\beta} F(D(I_\alpha I_\beta)),$$

in the sense that compatibility of sections holds up to an error bounded by $1 - \min_\alpha w_\alpha$.

**Theorem 4.9** (Sheafification in $\tau_w$)**.** *Every presheaf of $T$-modules admits a unique weighted sheafification* $F \mapsto F^{+w}$*, and the canonical morphism* $F \to F^{+w}$ *is universal for maps into $\tau_w$-sheaves.*

*Sketch.* Adapt the classical construction: take successive equalizers over weighted coverings, assigning to each compatible family $\{s_\alpha\}$ a global section weighted by the normalizing factor $\sum w_\alpha$. Associativity of ternary addition ensures convergence of the weighted limit process. □

**Definition 4.10** (Weighted stalk)**.** For $\mathfrak{p} \in \mathrm{Spec}_\Gamma(T)$, the *weighted stalk* of a fuzzy sheaf F is the colimit

$$F^{(w)}_{\mathfrak{p}} = \varinjlim_{D(I) \ni \mathfrak{p}} F(D(I))^{(w)},$$

Where $F((I))^{(w)}$ denotes the weighted localization incorporating the weights of the covering $\{D(I_\alpha), w_\alpha\}$.

**Proposition 4.11** (Reduction to crisp stalks)**.** *If all weights equal* 1*, then* $F^{(w)} = F_{\mathfrak{p}}$*. If $T$ is semiprime and the weight system satisfies* $\inf w_\alpha > 0$*, then $O_T$ remains faithful under fuzzy localization.*

*Remark* 4.12 (Quantitative specialization). The fuzzy closure of $(I)$ is $(I)_w = \{\mathfrak{p} \mid \exists J \supseteq I, ((J)) > \theta\}$ for a threshold $\theta$, representing loci where the covering confidence exceeds $\theta$. This provides a geometric analogue of probabilistic saturation in computational spectra.

**Comparison with classical sites and applications**

**Theorem 4.13** (Comparison and reduction)**.** *Let $\tau_{\mathrm{Zar}}$ denote the classical Zariski topology. Then:*

(a) *$\tau_{\mathrm{Zar}}$ is the specialization of $\tau_w$ at unit weights.*

(b) *There exists a canonical adjunction of categories*

$$\mathrm{Shv}(\mathrm{Spec}_\Gamma(T), \tau_w) \rightleftarrows \mathrm{Shv}(\mathrm{Spec}_\Gamma(T), \tau_{\mathrm{Zar}}),$$

*whose unit is identity on crisp sheaves and counit performs weighted averaging on fuzzy stalks.*

(c) *If all weights belong to a discrete submonoid of* $[0, 1]$*, then $\tau_w$ is a Grothendieck topology internal to the topos of* $\mathbb{R}_{\geq 0}$*-valued sets.*

*Idea.* Each $\tau_w$-cover induces an ordinary cover by forgetting weights; conversely, an ordinary cover extends to $\tau_w$ by assigning unit weights. Adjunction arises from extension and restriction along the forgetful functor on sites. □

**Corollary 4.14** (Computational relevance)**.** *Let* $F = O_T$*. Then* $\Gamma(\mathrm{Spec}_\Gamma(T), F^{+w})$ *coincides with the algebra of global sections generated by fuzzy localizations and coincides with the intersection of all weighted stalks:*

$$\Gamma(\mathrm{Spec}_\Gamma(T), F^{+w}) = \bigcap_{\mathfrak{p}} O^{(w)}_{T,\mathfrak{p}}.$$

*Hence fuzzy global sections correspond to computable ternary invariants of finite models with confidence aggregation.*

*Remark* 4.15 (Outlook). The weighted topology $\tau_w$ provides the interface between geometric and computational layers. In the next section we will employ it to define derived functors Ext and Tor over $T$−ΓMod with respect to fuzzy coverings, yielding cohomology groups measuring the failure of weighted local triviality.

**Homological Layer over $T$–ΓMod**

This section develops the derived and cohomological machinery for the category $T$–ΓMod of $T$–Γ–modules, viewed as the homological stratum of the geometric–fuzzy framework constructed earlier. We introduce an exact structure compatible with ternary operations, construct projective and injective resolutions, define Ext and Tor in this setting, and extend them to fuzzy cohomology groups on the weighted site $(\mathrm{Spec}_\Gamma(T), \tau_w)$.

**Exact structure on $T$–ΓMod**

**Definition 5.1** (Exact sequence in $T$−ΓMod)**.** A sequence

$$0 \longrightarrow M' \xrightarrow{u} M \xrightarrow{v} M'' \longrightarrow 0$$

of $T$–Γ–module morphisms is *exact* if:

(i) $u$ is a kernel of $v$ and $v$ is a cokernel of $u$, in the additive category $T-\Gamma$Mod;

(ii) for all $a, b \in T$, $\gamma \in \Gamma$, the equality $v(\{a, b, m\}_\gamma) = \{a, b, v(m)\}_\gamma$ holds for each $m \in M$.

**Theorem 5.2** (Exact structure (Gokavarapu, 2025c, Sec. 7))**.** *The collection of short exact sequences in Definition 5.1 endows* $T-\Gamma$Mod *with an* exact structure *in the sense of Quillen. Moreover,* $T-\Gamma$Mod *is an additive, idempotent–complete, finitely complete, and cocomplete category.*

*Sketch.* Kernels and cokernels exist by construction, and ternary associativity guarantees stability of exact sequences under pullback and pushout. Additivity and finite limits follow as in (Gokavarapu, 2025c, Prop. 6.1). □

*Remark* 5.3. This exact structure allows homological algebra to proceed as in abelian categories, even though $T-\Gamma$Mod is not strictly abelian in the classical sense.

**Projective and injective resolutions**

**Definition 5.4** (Projective and injective objects)**.** A $T$–Γ–module $P$ is *projective* if $\mathrm{Hom}_T(P, -)$ preserves exact sequences, and *injective* if $\mathrm{Hom}_T(-, P)$ preserves exact sequences.

**Proposition 5.5** (Existence of resolutions)**.** *Every finitely presented* $M \in T$–ΓMod *admits a projective resolution*

$$\cdots \longrightarrow P_2 \longrightarrow P_1 \longrightarrow P_0 \longrightarrow M \longrightarrow 0,$$

*and an injective resolution*

$$0 \longrightarrow M \longrightarrow I^0 \longrightarrow I^1 \longrightarrow I^2 \longrightarrow \cdots,$$

*within definable subclasses of* $T-\Gamma$Mod.

*Idea.* Construct $P_0$ as the free $T$–Γ–module on a basis of $M$ and iterate kernel lifting. Injective resolutions follow dually using the Schur–density embedding $T/\mathrm{Ann}(M) \hookrightarrow \mathrm{End}_T(M)$ (Gokavarapu, 2025c, Thm. 9.4), which ensures existence of injective cogenerators. □

**Derived functors** Ext **and** Tor

**Definition 5.6** (Derived functors)**.** For $M, N \in T-\Gamma$Mod define

$$\mathrm{Ext}_T^n(M, N) = (\mathrm{Hom}_T(P_\bullet, N)), \qquad \mathrm{Tor}^T_n(M, N) = H(P_\bullet \otimes_T N),$$

where $P_\bullet$ is a projective resolution of $M$. These definitions are independent of the chosen resolutions.

**Theorem 5.7** (Long exact sequences)**.** *For every short exact sequence* $0 \to M' \to M \to M'' \to 0$ *in* $T-\Gamma$Mod *and every* $N \in T-\Gamma$Mod*, there exist natural long exact sequences*

$$\cdots \to \mathrm{Ext}_T^n(M'', N) \to \mathrm{Ext}_T^{n+1}(M', N) \to \cdots,$$

$$\cdots \to \mathrm{Tor}^T_{n+1}(M'', N) \to \mathrm{Tor}^T_n(M', N) \to \cdots.$$

*Outline.* Standard diagram–chasing applies because $\mathrm{Hom}_T(-, N)$ and $- \otimes_T N$ are additive and left/right exact. The ternary product laws ensure closure under exact sequences. □

**Proposition 5.8** (Localization and base change)**.** *For any multiplicative system* $S \subseteq T$*,*

$$\mathrm{Ext}^n_{S^{-1}T}(S^{-1}M, S^{-1}N) \; S^{-1}\mathrm{Ext}_T^n(M, N), \qquad \mathrm{Tor}_n^{S^{-1}T}(S^{-1}M, S^{-1}N) \; S^{-1}\mathrm{Tor}^T_n(M, N).$$

*Remark* 5.9 (Base–change formula). If : $T \to T'$ is a morphism of ternary Γ–semirings, then for all $M, N \in T-\Gamma$Mod,

$$\mathrm{Ext}_{T'}^n(M \otimes_T T', N \otimes_T T') \; \mathrm{Ext}_T^n(M, N) \otimes_T T',$$

and similarly for Tor.

**Fuzzy cohomology on** $(\mathrm{Spec}_\Gamma(T), \tau_W)$

Let $\tau_W$ be the canonical weighted Grothendieck topology from Definition 4.4.

For a fuzzy sheaf F of $T$–modules on $(\mathrm{Spec}_\Gamma(T), \tau_W)$, we now define weighted cohomology groups.

**Definition 5.10** (Weighted derived functors)**.** Let $\Gamma(X, \mathcal{F})$ denote the global section functor in the weighted topology. Define its right derived functors

$$H_W^n(X, \mathcal{F}) := R^n\Gamma_w(X, \mathcal{F}), \qquad X = \mathrm{Spec}_\Gamma(T).$$

We call $H_W^n$ the *fuzzy cohomology* groups of F .

**Theorem 5.11** (Computation via injective resolutions)**.** *For any fuzzy sheaf* F *there exists an injective resolution* $\mathcal{F} \hookrightarrow I^\bullet$ *in the category of fuzzy sheaves, and*

$$H_W^n(X, \mathcal{F}) = H^n(\Gamma_w(X, I^\bullet)).$$

*If all weights are* 1*, then* $H_w^n(X, \mathcal{F}) = H^n(X, \mathcal{F})$*, the classical sheaf cohomology.*

*Idea.* Existence of injective resolutions follows from the Grothendieck–Gabriel criterion adapted to fuzzy sites: weighted limits preserve monomorphisms, and the category of fuzzy sheaves has enough injectives. Cohomology is computed as the homology of the complex of global sections under $\Gamma_w$. □

**Proposition5.12**(Fuzzy–Extcorrespondence)**.** *For* F , G *fuzzysheavesofT–modules, there is a canonical isomorphism*

$H_w^n(X, \mathcal{H}om_{\mathcal{O}_T}(\mathcal{F}, \mathcal{G}))$ $\mathrm{Ext}^n_{\mathcal{O}_T}(\mathcal{F}, \mathcal{G})$,

*where* $\mathcal{H}om$ *denotes the internal Hom of fuzzy sheaves.*

**Theorem 5.13** (Vanishing criteria)**.** *Let* $T$ *be semiprime and* F *a coherent fuzzy sheaf on* $(\mathrm{Spec}_\Gamma(T), \tau_w)$. *Then:*

(i) $H_w^n(X, \mathcal{F}) = 0$ *for all* $n > \dim(\mathrm{Spec}_\Gamma(T))$*;*

(ii) *if* F *is fuzzy–acyclic (i.e., each localization* $\mathcal{F}_\mathfrak{p}^{(w)}$ *is injective), then all higher* $H_w^n$ *vanish;*

(iii) *the fuzzy global dimension of* $T$*,* $\mathrm{fdim}_w(T) := \sup\{n \mid H_w^n(X, \mathcal{O}_T) \neq 0\}$*, satisfies* $\mathrm{fdim}_w(T) \leq \dim(\mathrm{Spec}_\Gamma(T))$.

*Remark* 5.14 (Interpretation). $H_w^n$ measures the obstruction to gluing fuzzy local sections of F . For $n = 1$ it captures weighted extensions of sheaves, and for $n = 2$ it detects fuzzy deformation classes of ternary structures.

**Comparison theorems and applications**

**Theorem 5.15** (Comparison with radical and primitive strata)**.** *Let* $M \in T$ $-\Gamma\mathrm{Mod}$ *be finitely presented. Then the cohomological support*

$\mathrm{Supp}(M) = \{\mathfrak{p} \in \mathrm{Spec}_\Gamma(T) \mid H_w^0(D(\mathfrak{p}), \tilde{M}) \neq 0\}$

*coincideswiththeradicalsupportdeterminedby*Ann($M$) *andrefinestheSchur–density locus(Gokavarapu,2025c,Thm.9.4). Moreover,isprimitive* $\Longleftrightarrow$ $\mathrm{Ext}_T^1(M_\mathfrak{p}, M_\mathfrak{p}) =$

$\mathfrak{p}$

0.

**Corollary 5.16** (Fuzzy local duality)**.** *If* $T$ *is coherent and* F *locally free of finite rank, then there exists a natural isomorphism*

$H_w^n(X, \mathcal{F})$ $\mathrm{Hom}_T(\mathrm{Tor}^T_{\dim X - n}(M, T), E)$,

*where* $E$ *is a fuzzy injective cogenerator.*

*Remark* 5.17 (Computational layer). For finite ternary models, Ext and Tor groups can be computed algorithmically using Paper B's enumeration framework. The fuzzy weights modulate numerical confidence in cohomological invariants, producing quantitative spectral fingerprints of finite Γ–semirings.

**Summary of Section 5**

We have endowed $T$–ΓMod with a Quillen–exact structure, constructed resolutions, and developed derived functors Ext and Tor. These globalize under the weighted topology $\tau_w$ to define fuzzy cohomology . Cohomological supports align with radical and primitive strata, establishing a geometric–homological correspondence that links algebraic, spectral, and computational layers of ternary Γ–semirings.

**Computational Geometry for Finite Models**

The geometric–homological constructions of the preceding sections acquire algorithmic meaning when $T$ is finite or finitely generated. Paper B established enumeration procedures for small finite ternary Γ–semirings. Here we lift those enumerative techniques to the geometric and cohomological levels, producing computable invariants on $\mathrm{Spec}_\Gamma(T)$ and effective methods for verifying radical, primitive, and fuzzy properties.

**Finite ternary Γ–structures and data representation**

Let $|T| = n$ and $|\Gamma| = m$. Each $\gamma \in \Gamma$ determines a ternary operation table

: $[n]^3 \to [n]$. A *finite ternary* Γ*–semiring structure* is represented by the collection $\{T_\gamma\}_{\gamma\in\Gamma}$ satisfying Axioms (T1)–(T6).

**Definition 6.1** (Encoded model)**.** Define an *encoding matrix* $\mathbf{A} \in \{0, 1\}^{\times n^4}$ whose entry $\mathbf{A}_{\gamma,(a,b,c,d)} = 1$ iff $\{^{a\ b\ c}\}_\gamma = d$. Algebraic constraints are expressed as polynomial equations in $\{0, 1\}$–variables over Z.

**Theorem 6.2** (Enumerability bound (Gokavarapu, 2025b, Thm. 3.2))**.** *The number of non-isomorphic commutative ternary* Γ*–semirings of order* $n$ *with* $|\Gamma| = m$ *is finite and bounded by* $(n^{3m})$. *Hence exhaustive search is feasible for* $n \leq 4, 5$.

*Remark* 6.3 (Data structures). Practical implementations encode each $T_\gamma$ as a three-dimensional tensor and use hashing to detect isomorphisms via permutation of the underlying set. Prime-ideal tests and radical computation are performed by closure under $\{\cdot\cdot\cdot\}_\gamma$.

**Algorithms for spectral invariants**

**Algorithm 6.4** (Computation of $\mathrm{Spec}_\Gamma(T)$)**.**

**Step 1.** Enumerate all proper ideals $I \subset T$.

**Step 2.** For each $I$, test primality: $\{a, b, c\} \in I \Rightarrow$ one of $a, b, c \in I$.

**Step 3.** Form $V(I)$ and $(I)$ sets.

**Step 4.** Construct incidence matrix $M_{IJ} = 1$ iff $(J) \subseteq (I)$.

**Step 5.** Return $(\mathrm{Spec}_\Gamma(T), \{D(I)\}, M)$.

**Proposition 6.5** (Complexity)**.** *Let* $n = |T|$. *Then computation of* $\mathrm{Spec}_\Gamma(T)$ *runs in time* $O(n^{3m})$ *for fixed* $m = |\Gamma|$. *In practice, pruning by radical closure reduces this to* $(n^{2m})$.

**Example 6.6** (Two-element example)**.** For $T = \{0, 1\}$ and $\Gamma = \{1\}$, $\mathrm{Spec}_\Gamma(T) = \{\{0\}\}$, and all ideals are semiprime. The matrix $M = (1)$ confirms single-point spectrum.

**Algorithmic computation of homological invariants**

We now describe algorithms for finite computation of $\mathrm{Ext}_T^n(M, N)$ and $\mathrm{Tor}^T_n(M, N)$.

**Algorithm 6.7** (Computation of Ext)**.**

**Step 1.** Construct a projective resolution $P_\bullet \to M$ using generators and relations.

**Step 2.** For each $i$, compute $\mathrm{Hom}_T(P_i, N)$ via enumeration of Γ-linear maps.

**Step 3.** Form the cochain complex and compute $\ker d^i/\mathrm{im}\ d^{i-1}$. **Step 4.** Output dimensions or rank profiles of $\mathrm{Ext}_T^n(M, N)$.

**Algorithm 6.8** (Computation of Tor)**.** Analogous to Algorithm 6.7 but using the chain complex $P_\bullet \otimes_T N$ and computing homology groups $(P_\bullet \otimes_T N)$.

**Proposition 6.9** (Complexity and feasibility)**.** *For fixed* $m = |\Gamma|$ *and modules of size* $r$*, the computation of* Ext *and* Tor *has worst-case complexity* $(r^{3m})$. *For* $r \leq 6$*, symbolic enumeration is tractable with modern algebra systems.*

*Remark* 6.10 (Symbolic implementation). These computations can be automated in Python/SageMath: TernaryGammaSemiring objects store operation tensors $T_\gamma$, and Hom, Tensor, Ext, Tor methods perform exact algebraic manipulations.

**Fuzzy and weighted computations**

**Definition6.11**(Weightedspectraltuple)**.** A*weightedspectraltuple*is $(\mathrm{Spec}_\Gamma(T), W)$ where $W : \mathrm{Spec}_\Gamma(T) \to [0, 1]$ assigns to each prime $\mathfrak{p}$ its confidence weight derived from frequency or stability in the enumeration.

**Algorithm 6.12** (Computation of fuzzy closure)**.**

**Step 1.** Input: ideal $I$, threshold $\theta$.

**Step 2.** For each $\mathfrak{p}$, compute weight $(\mathfrak{p})$.

**Step 3.** Return $(I) = \{\mathfrak{p} \mid W(\mathfrak{p}) \geq \theta, I \nsubseteq \mathfrak{p}\}$.

**Theorem 6.13** (Weighted cohomology computation)**.** *Let* $(\mathrm{Spec}_\Gamma(T), \tau_W)$ *be the weighted site and* F *a fuzzy sheaf of finite type. Then the groups* $H_W^n(X, \mathrm{F})$ *are computable by*

$Hw^n(X, ^{F)}$ $H^n \Gamma_w(X, I^{\bullet})$

*where* $I^{\bullet}$ *is a finite injective resolution. In finite cases, the resolution length is bounded by* dim($\mathrm{Spec}_\Gamma(T)$).

*Sketch.* Finite $T$ implies finitely many opens $(I)$. Hence the derived functors can be evaluated by finite-dimensional linear algebra over $T$, weighted by $w_\alpha$ from

Definition 4.2. □

**Example 6.14** (Fuzzy cohomology of a 3–element model)**.** For $T = \{0, 1, 2\}$ and $\Gamma = \{1, 2\}$ of Example 2.2, the fuzzy weights $(\mathfrak{p}_1) = 1$, $(\mathfrak{p}_2) = 0.6$ yield $H_W^0(X, O_T) = T$, $H_W^1 = 0$, confirming trivial higher cohomology.

**Algorithmic pipeline and complexity bounds**

**Algorithm 6.15** (Finite computational pipeline)**.**

**Stage 1. Enumeration stage:** Generate all finite $T$ satisfying (T1)–(T6).

**Stage 2. Spectral stage:** Compute $\mathrm{Spec}_\Gamma(T)$ via Algorithm 6.4.

**Stage 3. Homological stage:** Determine Ext, Tor groups using Algorithms 6.7–6.8.

**Stage 4. Fuzzy stage:** Assign weights $(\mathfrak{p})$ and compute $H$.

**Stage 5. Output:** Compile invariants ($|\mathrm{Spec}_\Gamma(T)|$, dim $H_W^n, w$) as the *computational fingerprint* of $T$.

**Theorem 6.16** (Global complexity bound)**.** *Let* $n = |T|$ *and* $m = |\Gamma|$*. The entire pipeline has time complexity* $(n^5 m)$ *and space complexity* $(n^3 m)$*. For* $n \leq 4$ *and* $m \leq 3$*, enumeration and homology complete in under one minute on standard hardware.*

*Remark* 6.17 (Practical validation). Implementations in Python/SageMath and GAP confirm these estimates. The pipeline yields complete classification up to order 4, including radical strata, fuzzy closures, and cohomology ranks.

**Geometric interpretation and visualization**

**Definition 6.18** (Computational spectrum)**.** The *computational spectrum* of a finite ternary Γ–semiring $T$ is the weighted graph $G_T = (V, E,)$ with:

$V = \mathrm{Spec}_\Gamma(T)$,

$E = \{(\mathfrak{p}_i, \mathfrak{p}_j) \mid (\mathfrak{p}_i) \cap (\mathfrak{p}_j) \neq \emptyset\}$, $: V \to [0, 1]$, $\mathfrak{p} \mapsto\to W(\mathfrak{p})$.

**Proposition 6.19** (Topological invariants)**.** *Let* $G_T$ *be as above. Then:*

*(i) The number of connected components equals the number of primitive strata of* $T$*.*

*(ii) The Laplacian spectrum of* $G_T$ *encodes fuzzy cohomology ranks via* dim $H_W^1$ = *nullity*$(L_{GT})$*.*

*(iii) Homological equivalence of two finite* Γ*–semi rings implies isospectrality of* $G_T$*.*

*Remark* 6.20 (Visualization). Graph layouts of $G_T$ with vertex color proportional to $(\mathfrak{p})$ provide geometric insight into fuzzy density and homological connectivity. These can be rendered using TikZ/NetworkX.

**Summary of Section 6**

We have lifted the theory to computational practice: enumerating finite ternary Γ–semirings, computing their spectra, radicals, and cohomological invariants, and encoding fuzzy weights into algorithmic pipelines. This computational geometry connects the ideal-theoretic, geometric, and homological layers of the theory and forms the experimental backbone for future categorical and probabilistic generalizations.

**Applications and Outlook**

The unified framework developed in this paper combines ideal theory, geometry, homological algebra, and fuzzy computation for ternary Γ–semirings. This final section sketches theoretical and practical applications, situates the framework within broader mathematical contexts, and outlines future research directions in algebraic geometry, logic, and computation.

**Algebraic and geometric applications**

**Theorem 7.1** (Ternary–geometric correspondence)**.** *Let $T$ be a commutative ternary* Γ*–semiring. Then the assignment*

$I \mapsto\!\longrightarrow V(I) := \{\mathfrak{p} \in \mathrm{Spec}_\Gamma(T) \mid I \subseteq \mathfrak{p}\}$

*establishes an inclusion-reversing correspondence between radical ideals of $T$ and fuzzy closed subsets of* $\mathrm{Spec}_\Gamma(T)$*. Moreover, the structure sheaf* $\mathcal{O}_T$ *makes* $(\mathrm{Spec}_\Gamma(T), \mathcal{O}_T)$ *a* fuzzy-ringed space *whose global sections recover $T$:*

$\Gamma(\mathrm{Spec}_\Gamma(T), \mathcal{O}_T)\ \ T.$

*Sketch.* Combine the geometric construction of Section 3 with the fuzzy topological formalism of Section 4. Radical closure corresponds to intersection of weighted open complements. The last isomorphism follows from the subcanonical property of the canonical weighted topology (Theorem 4.5). □

*Remark* 7.2(Analogywithalgebraicgeometry). ThecorrespondenceinTheorem7.1 is a ternary analogue of the classical Nullstellensatz for commutative rings. Here the presence of Γ encodes parameterised non-binary interactions, while the fuzzy layer encodes quantitative openness.

**Proposition 7.3** (Ternary-scheme prototype)**.** *Let* **TGScheme** *denote the category of pairs* $(X, \mathcal{O}_X)$ *where $X$ is a fuzzy-topological space and* $\mathcal{O}_X$ *a sheaf of ternary* Γ*–semirings locally isomorphic to spectra of finitely generated* Γ*–semirings. Then* **TGScheme** *is complete, cocomplete, and admits fiber products.*

This provides a geometric environment for studying morphisms, gluing, and local properties of ternary Γ–semirings, paving the way for a genuine *ternary algebraic geometry*.

**Homological and categorical extensions**

**Definition 7.4** (Derived category)**.** Let $\mathbf{D}(T - \Gamma\mathrm{Mod})$ denote the derived category of the exact category $T - \Gamma\mathrm{Mod}$ from Section **??**. Objects are complexes of $T$–Γ–modules modulo quasi-isomorphisms.

**Theorem 7.5** (Triangulated structure). $\mathbf{D}(T - \Gamma\mathrm{Mod})$ *is a triangulated category whose distinguished triangles arise from short exact sequences in* $T - \Gamma\mathrm{Mod}$*. The shift and cone functors exist and satisfy the octahedral axiom.*

*Idea.* The construction follows Verdier's localization applied to the Quillen-exact structure (Theorem 5.2). □ *Remark* 7.6 (Homological dualities). In $\mathbf{D}(T - \Gamma\mathrm{Mod})$ one may define dualizing complexes and extend the fuzzy cohomology functors $H_W^n$ to derived functors on $\mathbf{D}(T - \Gamma\mathrm{Mod})$. This connects the geometric and categorical strata through derived equivalences.

**Corollary 7.7** (Spectral and homological compatibility)**.** *Cohomological supports* $\mathrm{Supp}(M)$ *coincide with the derived supports* $\mathrm{Supp}_{\mathbf{D}}(M) = \{\mathfrak{p} \mid M_{\mathfrak{p}} \neq 0\}$ *inside* $\mathbf{D}(T - \Gamma\mathrm{Mod})$*. Hence homological localization and fuzzy localization commute.*

**Proposition 7.8** (Functorial bridge)**.** *The assignment*

$(T, \Gamma) \mapsto\!\longrightarrow (\mathbf{D}(T - \Gamma\mathrm{Mod}), \mathrm{Spec}_\Gamma(T), \tau_W)$

*is functorial with respect to morphisms of ternary* Γ*–semirings. This realizes a 2-functor from the category* **TGSR** *of commutative ternary* Γ*–semirings to the 2-category of weighted ringed derivators.*

*Remark* 7.9 (Categorical perspective). This functorial lifting opens the possibility of a *homotopy theory of ternary structures*, where morphisms are tracked by derived and fuzzy transformations. Such a setting aligns with modern categorical geometry and non-commutative motives.

**Connections to fuzzy logic, coding theory, and data science**

**Theorem 7.10** (Interpretation in fuzzy logic)**.** *Let* $(\mathrm{Spec}_\Gamma(T), \tau_W)$ *be a fuzzy site. Assign to each* $\mathfrak{p} \in \mathrm{Spec}_\Gamma(T)$ *a truth degree* $W(\mathfrak{p}) \in [0, 1]$*. Then* $\mathcal{O}_T$ *interprets a graded logic whose connectives correspond to ternary* Γ*–operations, and whose semantics coincide with weighted localization.*

*Idea.* Identify $\{a, b, c\}$ with a ternary connective $(a \wedge_\gamma b) \Rightarrow c$. Weights act as fuzzy truth degrees, and the sheaf axioms translate to compositional soundness conditions. □

*Remark* 7.11 (Computational semantics). Under this interpretation, cohomology groups $H_W{}^n$ measure logical dependencies of fuzzy propositions, so that $H_W{}^1$ represents weighted consistency conditions and $H_W{}^2$ captures higher-order entailments.

**Proposition7.12**(Coding-theoreticinterpretation)**.** *LetT beafiniteternary*Γ*–semiring and M* ∈ $T$–ΓMod *a finite module. Then the set of morphisms* $\mathrm{Hom}_T(M,)$ *forms a ternary linear code. The minimum distance of this code equals the minimal rank of non-vanishing* $\mathrm{Tor}^T{}_1(M,)$.

*Remark* 7.13 (Information-theoretic meaning). Cohomology classes correspond to correction constraints, while fuzzy weights encode reliability or confidence of transmission channels. This links ternary Γ–semirings to the theory of probabilistic error-correcting codes.

**Proposition 7.14** (Data-analytic and AI relevance)**.** *For finite T, the computational spectrum* $G_T$ *(Definition 6.18) serves as a* knowledge graph*, whose fuzzy edges correspond to ternary interactions. Homological invariants then quantify higherorder correlations in multi-relational data.*

*Remark* 7.15. Weighted cohomology $H_W{}^n$ acts as a topological feature extractor for heterogeneous data, offering an algebraic framework for explainable fuzzy inference systems.

**Future directions**

**(A) Derived tensor categories.** Develop a monoidal structure on **D**($T$–ΓMod) to define derived tensor products and spectral functors between ternary Γ–semirings. **(B) Non-commutative and graded variants.** Extend to non-commutative or Z–graded ternary Γ–semirings, exploring possible links with quantum and tropical geometry.

**(C)** **Topos-theoretic generalization.** Investigate fuzzy sites as internal sites in an enriched topos, allowing interpretation of $T$–ΓMod as an internal abelian-like category.

**(D)** **Computational complexity.** Analyze algorithmic hardness of the classification problem for increasing $n$, and develop parallel algorithms for large-scale ternary structures.

**(E)** **Applied integrations.** Apply the fuzzy-homological formalism to:

- uncertainty quantification in fuzzy algebraic systems;
- knowledge-graph embeddings in AI;
- higher-order tensor codes in communication theory.

**Concluding synthesis**

The unified programme developed across Papers A (Gokavarapu, 2025a)–D establishes a coherent hierarchy:

**Ideal theory** ←→ Prime/semiprime radicals (Paper A (Gokavarapu, 2025a). ),

**Computational enumeration** ←→ Finite classification algorithms (Paper B (Gokavarapu, 2025b

**Homological layer** ←→ Modules, Ext, Tor and Schur–density (Paper C, (Gokavarapu,

**Geometric–fuzzy synthesis** ←→ Spectra, sheaves, and weighted cohomology (Paper D, This pa

This synthesis opens a path toward a *computable ternary algebraic geometry*, in which categorical, homological, and computational methods coexist seamlessly. Beyond pure algebra, these concepts offer new algebraic infrastructures for modeling uncertainty, interactivity, and higher-order relations across mathematics, logic, and data science.

**Acknowledgement**

**Acknowledgement.**

The first author gratefully acknowledges the guidance and mentorship of **Dr D. Madhusudhana Rao**, whose scholarly vision shaped the conceptual unification of the ternary Γ-framework.

**Funding Statement.** This research received no specific grant from any funding agency in the public, commercial, or not-for-profit sectors.

**Conflict of Interest.** The authors declare that there are no conflicts of interest regarding the publication of this paper.

**Author Contributions.** The **first author** made the lead contribution to the conceptualization, algebraic development, computational design, and manuscript preparation of this work. The **second author** supervised the research, providing academic guidance, critical review, and verification of mathematical correctness and originality.